\documentclass[12pt]{article}

\usepackage{a4}
\usepackage[english]{babel}
\usepackage{graphicx}
\usepackage[centertags]{amsmath}
\usepackage{euscript}
\usepackage{amsfonts}
\usepackage{amssymb,amsbsy}
\usepackage{amsthm}
\usepackage{newlfont}

\newtheorem{thm}{Theorem}
\newtheorem{cor}[thm]{Corollary}
\newtheorem{lem}[thm]{Lemma}
\newtheorem{prop}[thm]{Proposition}
\newtheorem*{thm*}{Theorem}
\newtheorem{rem}{Remark}
\newtheorem*{rem*}{Remark}

\newtheorem{f}[thm]{Definition}

\newcommand{\ls}{\ldots}
\newcommand{\we}{\widetilde}
\newcommand{\st}{\subset}
\newcommand{\lb}{\big(}
\newcommand{\rb}{\big)}
\newcommand{\lfi}{\left\{ }
\newcommand{\rfi}{\right\} }

\newcommand{\xd}{\stackrel{\circ}{\b{x}}_{\c{D}^\iind{c}}}
\newcommand{\zd}{\stackrel{\circ}{\b{z}}_{\c{D}^\iind{c}}}
\newcommand{\xe}{\stackrel{\circ}{\b{x}}_{\c{E}^\iind{c}}}
\newcommand{\ze}{\stackrel{\circ}{\b{z}}_{\c{E}^\iind{c}}}

\newcommand{\bs}{\backslash}
\newcommand{\lla}{\longleftarrow}

\newcommand{\nn}{\noindent}

\def\b#1{\mathbf{\relax#1}}
\def\c#1{\mathcal{\relax#1}}
\def\bb#1{\mathbb{\relax#1}}

\def\rf#1{\mbox{\hspace{-1.5mm}\fontshape{ui}
            \selectfont{ \relax#1 }\hspace{-1.5mm}}}

\def\ind#1{{\mbox{\fontsize{9}{10.8}\selectfont$\relax#1 $} }}
\def\iind#1{{\mbox{\fontsize{8}{9.6}\selectfont$\relax#1 $} }}

\DeclareMathOperator*{\argmin}{argmin}

\begin{document}

\title{Efficient Determination of Gibbs  Estimators with
Submodular Energy Functions}

\author{Boris Zalesky\footnote{The research was supported by ISTC
B-517 grant.}}

\maketitle

\vspace{-0.3cm}

\begin{center}
United Institute of Information Problems\\
National Academy of Sciences\\
Surganov Street 6, Minsk, 220012, Belarus\\
e-mail address: zalesky@mpen.bas-net.by
\end{center}

\date{}

\vspace{0.5cm}

\begin{abstract}
Now the Gibbs estimators are  used to solve many problems in
computer vision,  optimal control, Bayesian statistics etc. Their
determination   is equivalent to  minimization of corresponding
energy functions. Therefore, the real use of the estimators is
restricted by opportunity to compute or evaluate satisfactory
values of the energy functions. New combinatorial techniques that
have been developed for the last years allowed efficient
employment of high dimensional Ising models with Boolean and
rational valued variables. In addition, graph cut methods were
proposed to minimize energy functions $ U=\sum g_\mu (x_\mu)+\sum
g_{\mu,\nu}(x_\mu,x_\nu)+ \sum
g_{\mu,\nu,\kappa}(x_\mu,x_\nu,x_\kappa) $ of Boolean variables
$x_\mu$ and $ U=\sum_\mu h(i_\mu)+\sum_{\mu,\nu}g(|i_\mu-j_\nu|) $
of integer variables $i_\mu$ with submodular functions $g$.

In the paper energy functions that often occur  in applications
(for instance, in estimating of images) are considered. They are
suppose to be depended on  variables taking values in  totally
ordered finite sets (finite sets of arbitrary numbers are among
them). The representation of functions by Boolean polynomials are
described. Fast methods of minimization of submodular Boolean
polynomials that correspond to the energy functions are developed.
For that the \emph{Modified algorithm of minimization of
submodular functions} (MSFM) is proposed. Also the \emph{minimum
graph cut technique} is presented for  submodular Boolean
polynomials that satisfy some additional condition. Among such
polynomials are, for instance, those with all coefficients before
nonlinear monomials negative as well as submodular Boolean
polynomials, which have all coefficients before monomials of
order more than two  positive.
\end{abstract}

\vspace{0.5cm}

\nn

\text{}\\
\nn
\newpage

\section{Introduction}
Nowadays Gibbs models are successfully employed in mathematical
statistics, statistical physics, image processing etc. Often these
models allow to get satisfactory solutions of various problems but
usually they turn out very hard to be solved efficiently. Finding
of Gibbs estimates is equivalent to the problem of identification
of minima of corresponding energy functions. In this paper we
present fast methods of finding  the Gibbs estimates for some
models that occur in application (for instance, in image
processing).

The Simulated annealing \cite{GG84,Bum01_1,Bum01_2,Z00}
was one of the first method
that was used to find the Gibbs estimates. In spite of its
slowness, it  remains a very  popular tool to evaluate the
estimates as well as other problems of general global optimization.
To speed up computation of the estimator of the Ising model Grieg,
Porteous and Seheult \cite{GPS89} in 1989 described the network flow
technique that allowed the exact Boolean minimizing  the Ising
energy function. In the last few years  several new approaches,
which enable efficient estimation \cite{Boyk01} or even efficient
minimizing some wide classes of functions of many Boolean or
integer variables, have been developed \cite{Ishik01,Kol01, Z02}.
These  approaches, which are based on the graph cut technique,
demonstrate possibility to solve real practical problems including
handling of 3D data for acceptable time  (for instance,
segmentation of 3D $200\times 200\times 300$ gray-scale image by
the such type method took 14min. of Pentiym 4 \ 1.6GH  CPU).

Recently \emph{the combinatorial algorithm of minimization of
submodular functions} (SFM) have been developed \cite{Iw02}. To
minimize a submodular function of $n$ variables it needs either
$O(n^5\log M)$ (where $M$ is some number depending on the function
considered) or $O(n^8)$ operations. Therefore, the SFM can not be
immediately used for finding of high dimensional Gibbs estimates
for its comparatively high computational expenses. We present a
modification of the SFM that can be useful for minimization of
submodular Boolean representations of energy functions, which
occur in applications, for instance, in image processing.

In order to use the SFM algorithm and the graph cut technique for
energy functions, which depend on variables  taking values in
finite totally ordered sets, we represent the functions by
Boolean polynomials. Then we develop methods that allow efficient
minimization of submodular Boolean polynomials of special
structure. Some Boolean polynomials admit finding parts of
Boolean vectors, which minimize these polynomials, by use only
parts of their monomials. We  illustrate this property with image
processing terminology. It means colours of far points of images
estimates are independent, i.e. a change of colours of some image
pixels in the estimate does not affect  colours of other far
pixels. Therefore, \emph{the modified SFM} (MSFM) operates, first,
not with an original Boolean polynomial but with separated its
parts, trying to find part of coordinates of a solution of the
original problem. Then  the values found are set in the original
polynomial and the procedure are repeated either with other parts
of the residual polynomial or with the entire residual polynomial.
At worst, it takes as many operations as the usual SFM and, at
best, the number of operations is even $O(n)$, where $n$ is a number
of Boolean variables.

As we mentioned above, new graph cut methods \cite{Kol01,Z02} are
now an efficient tool to determine Gibbs estimators with energy
functions $ U=\sum g_\mu (x_\mu)+\sum g_{\mu,\nu}(x_\mu,x_\nu)+
\sum g_{\mu,\nu,\kappa}(x_\mu,x_\nu,x_\kappa) $ that are Boolean
polynomials of the third order. We failed to  represent an
arbitrary submodular Boolean polynomial by a graph but we made it
for submodular Boolean polynomials of the special structure that
can be met in applications. The set of those polynomials
contains, for instance,
Boolean polynomials with all nonlinear monomials negative (of
course, such polynomials are submodular) and  submodular Boolean
polynomials with monomials of order three and higher positive.

\section{Representation of Functions by Boolean Polynomials}\label{s:rf}

Below the function $U(\b{x}),\ \b{x}=(x_1,\ldots,x_n)$ is supposed
to be dependent on variables that take values in an arbitrary
finite totally ordered set $\c{R}=\{r_0,r_1,\ldots,r_k\},$ $\
r_0\le r_1 \le \ldots\le r_k$. The multiindex
$\b{j}=(j_1,j_2,\ldots,j_n)$ is with coordinates $j_l\in
\{0,1,\ldots,k\}={\mathbb N}_k$. The  function of integer variables
$$
V(\b{j})=U(r_\ind{j_\iind{1}},r_\ind{j_\iind{2}},\ldots,
r_\ind{j_\iind{n}})
$$
is considered instead of the initial $U(\b{x})$.

To exploit the MSF and graph cut technique we, first, represent
the integer variables $j_\ind{i}$ as the sum of ordered Boolean ones
\begin{equation}\label{e:xdec}
j_\ind{i}=\sum_{l=1}^k x_\ind{i}(l),
\end{equation}
where $x_i(l)\in\{0,1\},\ x_i(1)\ge x_i(2)\ge\ldots\ge x_i(k).$ By
other words, we establish  one-to-one correspondence between each
variable $j_i$ and the ordered sequence of Boolean variables
$x_i(1)\ge x_i(2)\ge\ldots\ge x_i(k)$. After we represent the
function $V(\b{j})$ by a Boolean polynomial and then minimize
the polynomial.

Denote the Boolean vector $\b{x}(l)=(x_1(l),x_2(l),\ldots,
x_n(l))$, and for two vectors $\b{x}, \b{z}$ say $\b{x}\ge \b{z}$,
if  $x_i\ge z_i,\ i=1\div n$ (respectively, $\b{x}\ngeq\b{z}$, if
there is a couple  of indexes $i,j$ such that  $x_i < z_i$ and
$x_j > z_j$). Then
$$
V(\b{j})= V\left(\sum_{l=1}^k \b{x}(l)\right), \quad \b{x}(1)
\ge\b{x}(2)\ldots\ge\b{x}(k).
$$
For simplicity  let us start with the function of two variables
$V(j_1,j_2)$. It is not hard to check the equality
\begin{equation}\label{e:v1}
V(j_1,j_2)=V(0,j_2)+\sum_{\mu=1}^k(V(\mu,j_2)-V(\mu-1,j_2))x_1(\mu),
\end{equation}
where ordered Boolean variables $x_1(\mu)$ satisfy the equality
(\ref{e:xdec}). Now the expansion
\begin{multline}
V(j_1,j_2)\\
=V(0,0)+\sum_{\mu=1}^k(V(\mu,0)-V(\mu-1,0))x_1(\mu)
+\sum_{\nu=1}^k(V(0,\nu)-V(0,\nu-1))x_2(\nu)\\
+\sum_{\mu,\nu=1}^k
(V(\mu,\nu)-V(\mu-1,\nu)-V(\mu,\nu-1)+V(\mu-1,\nu-1))x_1(\mu)x_2(\nu),
\end{multline}
can be constructed by use of equality (\ref{e:v1})  for the second
variable $j_2$ of each difference $ V(\mu,j_2)-V(\mu-1,j_2). $

In general case, for $\b{j}\in \bb{N}_n$ we say
$$\Delta_{\ind{i}} V(\b{j})=V(\b{j})-V(\b{j}-\b{e}_\ind{i}),\
\b{e}_\ind{i}=(\overbrace{\underbrace{0,\ldots,0}_{i-1},1,0,\ldots,0}^n)$$
is difference derivative of the first order, and define the mixed
difference derivative of higher orders by the recursion
$$\Delta_{\ind{j}_\iind{1},\ind{j}_\iind{2},\ldots,
\ind{j}_\iind{l+1}}V(\b{j})=
\Delta_{\ind{j}_\iind{1},\ind{j}_\iind{2},\ldots,
\ind{j}_\ind{l}}V(\b{j})-\Delta_{\ind{j}_\iind{1},\ind{j}_\iind{2},
\ldots,\ind{j}_\iind{l}}V(\b{j}-\b{e}_{\ind{j}_\iind{l+1}}).
$$
Note that for any  permutation of indices
$\pi(j_1,j_2,\ldots,j_n)$
$$
\Delta_{\ind{\pi}(\ind{j}_\iind{1},\ind{j}_\iind{2},\ldots,
\ind{j}_\iind{l+1})}=\Delta_{\ind{j}_\iind{1},\ind{j}_\iind{2},\ldots,
\ind{j}_\iind{l+1}}.
$$
By analogy with 2D case
\begin{equation}\label{e:d1}
V(\b{j})=U(0,j_2,\ldots,j_n)+\sum_{\ind{\mu}=1}^k
\Delta_{1}V(\mu,j_2,\ldots,j_n))x_1(\mu).
\end{equation}
Let the $G_m$ be the projection operator onto the space of last
$n-m$ coordinate, i.e.
$$G_m\b{j}=(0,\ldots,0,j_{m+1},j_{m+2},\ldots,j_{n}).$$
Applying  formula (\ref{e:d1}) recursively $m$ times to each its
item  we get the equality
\begin{multline*}
V(\b{j})=V(G_m\b{j}) +\sum_{1\le \iind{l}\le m}\sum_{\mu=1}^k
\Delta_{l}V(G_m\b{j}+\mu\b{e}_l)x_l(\mu)\\
+\sum_{1\le \iind{l}_\iind{1}<\iind{l}_\iind{2}\le m}
\sum_{\ind{\mu}_{\iind{1}},\ind{\mu}_{\iind{2}}=1}^k
\Delta_{\iind{l}_\iind{1},\iind{l}_\iind{2}}
V(G_m\b{j}+\mu_{1}\b{e}_{l_1}+\mu_{2}\b{e}_{l_2})
x_{l_\ind{1}}(\mu_{\ind{1}})x_{l_\ind{2}}(\mu_{\ind{2}})+\ls\\
+\sum_ {\ind{\mu}_\iind{1},\ls,\ind{\mu}_\iind{m}=1}^k
\Delta_{1,2,\ls,m}
V(G_m\b{j}+\mu_{1}\b{e}_{1}+\mu_{2}\b{e}_{2}+\ls+\mu_{m}\b{e}_{m})
x_{{1}}(\mu_{{1}})x_{{2}}(\mu_{{2}})\ls x_{{m}}(\mu_{{m}})
\end{multline*}
Therefore, the following statement is valid.
\begin{prop}\label{e:pol}
The expansion of $V(\b{j)}$ into a  polynomial in ordered Boolean
variables $\b{x}(1)\ge\b{x}(2)\ge\ls\ge \b{x}(k)$ is of the form
\begin{multline}\label{e:ex}
V(\b{j})=\we{P}_V\left(\b{x}(1),\b{x}(2)\ldots,\b{x}(k)\right)\\
=V(\b{0}) +\sum_{m=1}^n\ \ \sum_\ind{1\le {l}_{1}<\ls<{l}_{m}\le
n}\ \ \sum_{\ind{\mu}_{\iind{1}},\ls,\ind{\mu}_{\iind{m}}=1}^k
\Delta_{{l}_{1},\ls,{l}_{m}} V\left(\sum_{\varkappa=1}^m
\mu_{\ind{\varkappa}}\b{e}_{l_\iind{\varkappa}}\right)
\prod_{\varkappa=1}^m
x_{l_\iind{\varkappa}}(\mu_{\ind{\varkappa}}).
\end{multline}
\end{prop}
The expansion (\ref{e:ex}) does not allow to turn immediately  to
the problem of the Boolean minimization because of ordered
variables. To avoid this disadvantage we will use another
polynomial
$$
P_V\left(\b{x}(1),\b{x}(2)\ldots,\b{x}(k)\right)=
\we{P}_V\left(\b{x}(1),\b{x}(2)\ldots,\b{x}(k)\right)
+\c{C}\sum_{\mu=1}^n\sum_{l=2}^k (x_\mu(l)-x_\mu(l-1))x_\mu(l),
$$
for sufficiently large constant $\c{C}>0$. This polynomial satisfy
the following properties.
\begin{prop}
The inequality  $\we{P}_V\le P_V$ holds true. Any collection of
Boolean vectors $\b{u^*}=(\b{x^*}(1),\b{x^*}(2)$
$\ldots,\b{x^*}(k))$ that minimizes $P_V$
$$
\b{u}^*=\argmin_{\b{x}(1),\b{x}(2),\ls,\b{x}(k)}
P_V\left(\b{x}(1),\b{x}(2)\ldots,\b{x}(k)\right)
$$
turns to be ordered (non-increasing). This collection minimizes
the polynomial $\we{P}_V$, and, therefore, the integer vector
$\b{j}^*=\sum_{l=1}^k \b{x^*}(l)$ minimizes the original function
$V(\b{j})$.
\end{prop}
Therefore, instead of minimization of $V(\b{j})$ we can consider the
problem of Boolean minimization of the polynomial
$$
P_V(\b{u})= a_0 + \sum_{m=1}^n\ \ \sum_\ind{1\le
{l}_{1}<\ls<{l}_{m}\le n}\ \
a_{l_\iind{1},l_\iind{2},\ls,l_\iind{m}} \prod_{i=1}^m
u_{l_\iind{i}},\qquad \dim(\b{u})={kn}.
$$
Because of the representation of functions $V(\b{j})$ by Boolean
polynomials developed we consider the problem of minimization of
general Boolean polynomials $P(\b{x})$.

\section{Submodular Boolean Polynomials}

For two Boolean vectors $\b{x,y}$ let us denote by
$\b{x}\vee\b{y}$ (respectively, by $\b{x}\wedge\b{y}$) the vector
with coordinates  $\max(x_i,y_i)$ (respectively, with coordinates
$\min(x_i,y_i)$).
\begin{f}\label{d:sub}
A function $f$ is called \emph{submodular} if for any Boolean
vectors $\b{x,y}$ it satisfies the inequality
$$
f(\b{x}\wedge\b{y})+f(\b{x}\vee\b{y})\le f(\b{x})+f(\b{y}).
$$
\end{f}

Let the set $\c{V}=\{1,\ls,n\}$, and for any $\c{D}\subset\c{V}$
the set $\c{D}^\ind{c}=\c{V}\backslash\c{D}$. Let the vector
$\b{x}_\c{D}$ with coordinated $x_i,\ i\in\c{D}$ be the restriction
of  $\b{x}$ to the set $\c{D}$, the abbreviation
$P(\b{x}_\c{D})=P(\b{x}_\c{D},\b{0}_{\c{D}^\ind{c}})$.

Below we  always suppose $P(\b{0})=0$. Therefore,
$P$ is of the form
\begin{equation}\label{e:p0}
P(\b{x})=  \sum_{m=1}^n\ \ \sum_\ind{1\le
{l}_{1}<\ls<{l}_{m}\le n}\ \
a_{l_\iind{1},l_\iind{2},\ls,l_\iind{m}} \prod_{i=1}^m
x_{l_\iind{i}}.
\end{equation}
\label{pij}

For any
$i,j\in\c{V}$ a polynomial $P(\b{x})$ can be written in the form
\begin{equation}\label{e:pij}
P(\b{x})= P(\b{x}_{\{i,j\}^\iind{c}}) +
x_iP_{i}(\b{x}_{\{i,j\}^\iind{c}})+
x_jP_{j}(\b{x}_{\{i,j\}^\iind{c}})+
x_ix_jP_{i,j}(\b{x}_{\{i,j\}^\iind{c}})
\end{equation}
 or in the form
$$
P(\b{x})=Q(\b{x})+L(\b{x}),
$$
where $Q(\b{x})$ is the polynomial consisting of monomials of
degree $2$ and higher, and $L(\b{x})$ is a linear polynomial. It
also can be represented as
$$
P(\b{x}) = P\{\c{D}\}(\b{x}) +
P(\b{x}_{\c{D}^\iind{c}}),\quad\c{D}\st \c{V},
$$
where each monomial of $P\{\c{D}\}(\b{x})$ contains at least one
variable $x_i,\  i\in\c{D}$ (and, in general, it depends on
$\b{x}_{\c{D}^\iind{c}}$), and $P(\b{x}_{\c{D}^\iind{c}})$ is
independent of $\b{x}_\c{D}$. The following useful properties can
be easily veryfied.
\begin{prop}\label{p:prop}
a) A Boolean polynomial $P(\b{x})$ is submodular if and only if
for any  vector $\b{x}$ and any couple of indices $i,j\in\c{N}$
the polynomials $P_{i,j}(\b{x}_{\{i,j\}^\iind{c}})\le 0$.

b) In order to  Boolean polynomial $P(\b{x})$ be submodular it is
necessary and sufficient its nonlinear part $Q(\b{x})$ will be
submodular.

c) If $P(\b{x})$ is submodular, then the polynomial $Q(\b{x})$ is
nonincreasing, in the sense that $Q(\b{y})\le Q(\b{x})$ for any
ordered pair of vectors $\b{y}\ge\b{x}$. Therefore,
$\argmin_{\b{x}}Q(\b{x})=\b{1}_{\c{V}}$.

d) If $P(\b{x})$ is submodular, then $P\{\c{D}\}(\b{x})$ is
submodular polynomial in $\b{x}_\c{D}$ for each fixed vector
$\b{x}_{\c{D}^\iind{c}}$.
\end{prop}

The sentences a), b) and d) are deduced directly from Definition
\ref{d:sub}. To prove c) let us  consider a couple of vectors:
$\b{x}$ with $x_i=0$ and $\b{y}=\b{x}+\b{e}_i$ and inequality
$$
Q(\b{0}_{\c{V}}) + Q(\b{y})\le Q(\b{e}_i) + Q(\b{x})
$$
that follows from Definition \ref{d:sub}. Since
$Q(\b{0}_{\c{V}})=Q(\b{e}_i)=0$ we have $Q(\b{y})\le Q(\b{x})$. To
prove c) in general case it is enough to repeat the last inequality
successively.

We need two properties of polynomials $P\{\c{D}\}(\b{x})$. The
first one is rather clear.
\begin{prop}\label{p:pd}
If $\b{x}'_\c{D}$ minimizes the polynomial $P\{\c{D}\}(\b{x})$ for
some fixed $\stackrel{\circ}{\b{x}}_{\c{D}^\iind{c}}$, then for
any set $\c{E}\st \c{D}$ the restriction $\b{x}'_\c{E}$ minimizes
the polynomial $ P\{\c{E}\}(\b{x}_\c{E},\b{x}'_{\c{D} \backslash
\c{E}},\stackrel{\circ}{\b{x}}_{\c{D}^\iind{c}}) $ with respect to
$\b{x}_\c{E}$, and vice versa, if the vector
$\widehat{\b{x}}_\c{E}$ minimizes $
P\{\c{E}\}(\b{x}_\c{E},\b{x}'_{\c{D} \backslash
\c{E}},\stackrel{\circ}{\b{x}}_{\c{D}^\iind{c}}), $ then the vector
$ (\widehat{\b{x}}_\c{E},\b{x}'_{\c{D} \backslash \c{E}}) $
minimizes $
P\{\c{D}\}(\b{x}_\c{D},\stackrel{\circ}{\b{x}}_{\c{D}^\iind{c}}). $
with respect to $\b{x}_\c{D}$.
\end{prop}

The second property concerns  monotonic (in some sense) dependence
of the solution $ \b{x}'_\c{D}=\argmin_{\ind{\b{x}}_\c{D}}
P\{\c{D}\}(\b{x}_\c{D},\stackrel{\circ}{\b{x}}_{\c{D}^\iind{c}}) $
on boundary vector $\stackrel{\circ}{\b{x}}_{\c{D}^\iind{c}}$ for
a submodular polynomial $P$. Let
$$
\c{M}({\stackrel{\circ}{\b{x}}_{\c{D}^\iind{c}}})= \left\{
\argmin_{\ind{\b{x}}_\c{D}}
P\{\c{D}\}(\b{x}_\c{D},\stackrel{\circ}{\b{x}}_{\c{D}^\iind{c}})
\right\}
$$
be the set of  minima of the polynomial $P\{\c{D}\}$ for fixed
$\stackrel{\circ}{\b{x}}_{\c{D}^\iind{c}}$. It can be shown that
for two ordered boundary vectors $\xd\le\zd$ can exist either
unordered $\b{x}'_\c{D}\nleq\b{z}'_\c{D}$ or non-increasing
$\b{x}'_\c{D}\ge\b{z}'_\c{D}$ solutions, but for any
$\b{x}'_\c{D}$ there exists $\b{z}'_\c{D}$ such that
$\b{x}'_\c{D}\leq\b{z}'_\c{D}$, and vice versa, for any
$\b{z}'_\c{D}$ there is $\b{x}'_\c{D}$ satisfying the same
inequality.
\begin{thm}\label{t:main}
Let a polynomial $P$ be submodular. Let for $\xd\le\zd$ solutions $
\b{x}'\in\c{M}({\stackrel{\circ}{\b{x}}_{\c{D}^\iind{c}}}) $ and $
\b{z}'\in\c{M}({\stackrel{\circ}{\b{z}}_{\c{D}^\iind{c}}}) $ do
not satisfy the inequality $\b{x}'_\c{D}\le\b{z}'_\c{D}$.

Then for the maximal nonempty set of the inverse order $\c{E}=\{i\in
\c{D}\ |\ x'_i=1,\ z'_i=0 \}$, such that
$\b{1}_\c{E}=\b{x}_\c{E}'>\b{z}'_\c{E}=\b{0}_\c{E}$, the vector $
(\b{z}'_\c{E},\b{x}'_{\c{D}\bs\c{E}} )=
(\b{0}_\c{E},\b{x}'_{\c{D}\bs\c{E}} )\in
\c{M}({\stackrel{\circ}{\b{x}}_{\c{D}^\iind{c}}}) $ as well as $
(\b{x}'_\c{E},\b{z}'_{\c{D}\bs\c{E}} )=
(\b{1}_\c{E},\b{z}'_{\c{D}\bs\c{E}} )\in
\c{M}({\stackrel{\circ}{\b{z}}_{\c{D}^\iind{c}}}) $
and, therefore,
$
(\b{z}'_\c{E},\b{x}'_{\c{D}\bs\c{E}} )\le
(\b{x}'_\c{E},\b{z}'_{\c{D}\bs\c{E}} ).
$
\end{thm}
\begin{proof}
For boundary vectors $\xd\le\zd$ suppose the condition
$\b{x}'_\c{D}\le\b{z}'_\c{D}$ is not valid, i.e. either
$\b{x}'_\c{D}\ge\b{z}'_\c{D}$ or $\b{x}'_\c{D}\nleq\b{z}'_\c{D}$.
In this case  the set $\c{E}=\{i\in \c{D}\ |\ x'_i=1,\ z'_i=0 \}$
is nonempty. Let us consider the chain of relations. The inequality
\begin{equation}\label{e:chain1}
P\{\c{E}\}(\b{z}'_\c{E},\b{z}'_{\c{D} \backslash
\c{E}},\stackrel{\circ}{\b{z}}_{\c{D}^\iind{c}}) \le
P\{\c{E}\}(\b{x}'_\c{E},\b{z}'_{\c{D} \backslash
\c{E}},\stackrel{\circ}{\b{z}}_{\c{D}^\iind{c}})
\end{equation}
holds true since the vector $\b{z}'_\c{E}$ minimizes $
P\{\c{E}\}(\ \cdot\ ,\b{z}'_{\c{D} \backslash
\c{E}},\stackrel{\circ}{\b{z}}_{\c{D}^\iind{c}}) $ (see,
Proposition \ref{p:pd}). The inequality
\begin{equation}\label{e:chain2}
P\{\c{E}\}(\b{x}'_\c{E},\b{z}'_{\c{D} \backslash
\c{E}},\stackrel{\circ}{\b{z}}_{\c{D}^\iind{c}}) \le
P\{\c{E}\}(\b{x}'_\c{E},\b{x}'_{\c{D} \backslash
\c{E}},\stackrel{\circ}{\b{x}}_{\c{D}^\iind{c}})
\end{equation}
will be proved later. The relation
\begin{equation}\label{e:chain3}
P\{\c{E}\}(\b{x}'_\c{E},\b{x}'_{\c{D} \backslash
\c{E}},\stackrel{\circ}{\b{x}}_{\c{D}^\iind{c}}) \le
P\{\c{E}\}(\b{z}'_\c{E},\b{x}'_{\c{D} \backslash
\c{E}},\stackrel{\circ}{\b{x}}_{\c{D}^\iind{c}})
\end{equation}
is verified the same way as (\ref{e:chain1}). Equalities
\begin{equation}\label{e:chain4}
P\{\c{E}\}(\b{z}'_\c{E},\b{x}'_{\c{D} \backslash
\c{E}},\stackrel{\circ}{\b{x}}_{\c{D}^\iind{c}}) =
P\{\c{E}\}(\b{z}'_\c{E},\b{z}'_{\c{D} \backslash
\c{E}},\stackrel{\circ}{\b{z}}_{\c{D}^\iind{c}}) =0
\end{equation}
are valid because $\b{z}'_\c{E}=\b{0}_\c{E}$. Therefore, relations
(\ref{e:chain1})-(\ref{e:chain4}) are, actually, equalities, and
the claim of the Theorem can be derived from Proposition
\ref{p:pd}.

Thus, the Theorem will be proved if we ground the inequality
(\ref{e:chain2}). For that we need additional notations.
One-element sets $\{j\}\in\c{V}$
 will be denoted by $j$. Instead of
$\c{E}\cup\{j\}$ and $\c{E}\bs\{j\}$ we will use  abbreviations
$\c{E}+j$ and $\c{E}-j$. For arbitrary sets $\c{A,B,C}\in\c{V}$ we
denote by $P\{\c{A},\c{B}|\c{C}\}(\b{x})$ polynomials that
consist  of monomials of $P(\b{x})$ with variables $x_i,\ i\in\c{A}
,\ x_j,\ j\in\c{B}$ but without $x_k,\ k\in\c{C}$.

First, consider Boolean vectors that for some $j\in\c{E}^\ind{c}$
satisfy the equation $\xe+\b{e}_j=\ze$ and the polynomial $
P\{\c{E}+ j\}(\b{x}_{\c{E}+j},
\stackrel{\circ}{\b{x}}_{\c{E}^\ind{c}- j}), $ which is a
submodular function of $\b{x}_{\c{E}+j}$. Let us write the
condition of submodularity of the polynomial for vectors
$\b{a}_{\c{E}+j}=(\b{0}_{\c{E}},1_j)$ and
$\b{b}_{\c{E}+j}=(\b{x}_{\c{E}},0_j)$ using the representation
\begin{multline*}
P\{\c{E}+ j\}(\b{x}_{\c{E}+j},
\stackrel{\circ}{\b{x}}_{\c{E}^\ind{c}- j})=
P\{\c{E}|\c{E}^\ind{c}\}(\b{x}_{\c{E}})+
P\{\c{E},j\}(\b{x}_{\c{E}+j},
\stackrel{\circ}{\b{x}}_{\c{E}^\ind{c}- j})+\\
P\{\c{E},\c{E}^\ind{c}-j|j\}(\b{x}_{\c{E}},
\stackrel{\circ}{\b{x}}_{\c{E}^\ind{c}- j})+
P\{j,\c{E}^\ind{c}-j|\c{E}\}(x_j,
\stackrel{\circ}{\b{x}}_{\c{E}^\ind{c}- j})+a_jx_j.
\end{multline*}
Since $ \b{a}_{\c{E}+j}\wedge \b{b}_{\c{E}+j}=\b{0}_{\c{E}+j} $ and
$ \b{a}_{\c{E}+j}\vee \b{b}_{\c{E}+j}=(\b{x}_{\c{E}},1_j), $ it
looks as
\begin{multline*}
P\{\c{E}|\c{E}^\ind{c}\}(\b{x}_{\c{E}})+
P\{\c{E},j\}(\b{x}_{\c{E}},1_j,
\stackrel{\circ}{\b{x}}_{\c{E}^\ind{c}- j})+\\
P\{\c{E},\c{E}^\ind{c}-j|j\}(\b{x}_{\c{E}},
\stackrel{\circ}{\b{x}}_{\c{E}^\ind{c}- j})+
P\{j,\c{E}^\ind{c}-j|\c{E}\}(1_j,
\stackrel{\circ}{\b{x}}_{\c{E}^\ind{c}- j})+a_j\le\\
\shoveleft{ P\{j,\c{E}^\ind{c}-j|\c{E}\}(1_j,
\stackrel{\circ}{\b{x}}_{\c{E}^\ind{c}- j})+a_j+}
\\
P\{\c{E}|\c{E}^\ind{c}\}(\b{x}_{\c{E}})+
P\{\c{E},j\}(\b{x}_{\c{E}},0_j,
\stackrel{\circ}{\b{x}}_{\c{E}^\ind{c}- j})+
P\{\c{E},\c{E}^\ind{c}-j|j\}(\b{x}_{\c{E}},
\stackrel{\circ}{\b{x}}_{\c{E}^\ind{c}- j}).
\end{multline*}
Therefore,
\begin{equation}\label{e:sube}
P\{\c{E},j\}(\b{x}_{\c{E}},1_j,
\stackrel{\circ}{\b{x}}_{\c{E}^\ind{c}- j})\le
P\{\c{E},j\}(\b{x}_{\c{E}},0_j,
\stackrel{\circ}{\b{x}}_{\c{E}^\ind{c}- j})=0.
\end{equation}
But
\begin{equation}\label{e:pe}
P\{\c{E}\}(\b{x}_{\c{E}}, \stackrel{\circ}{\b{x}}_{\c{E}^\ind{c}})=
P\{\c{E}|\c{E}^\ind{c}\}(\b{x}_{\c{E}})+
P\{\c{E},j\}(\b{x}_{\c{E}},x_j,
\stackrel{\circ}{\b{x}}_{\c{E}^\ind{c}- j})+
P\{\c{E},\c{E}^\ind{c}-j|j\}(\b{x}_{\c{E}},
\stackrel{\circ}{\b{x}}_{\c{E}^\ind{c}- j})
\end{equation}
Inequality (\ref{e:sube}), equation (\ref{e:pe}) and condition $
\stackrel{\circ}{\b{x}}_{\c{E}^\ind{c}- j}=
\stackrel{\circ}{\b{z}}_{\c{E}^\ind{c}- j} $ allows to conclude
\begin{equation}\label{e:pleq}
P\{\c{E}\}(\b{x}_{\c{E}},
\stackrel{\circ}{\b{x}}_{\c{E}^\ind{c}})\ge
P\{\c{E}\}(\b{x}_{\c{E}}, \stackrel{\circ}{\b{z}}_{\c{E}^\ind{c}}).
\end{equation}
To prove (\ref{e:pleq}) for an arbitrary pair $
\stackrel{\circ}{\b{x}}_{\c{E}^\ind{c}}\le
\stackrel{\circ}{\b{z}}_{\c{E}^\ind{c}} $ it is enough to use
monotone increasing sequence of boundary vectors $
\stackrel{\circ}{\b{x}}_{\c{E}^\ind{c}}<
\stackrel{\circ}{\b{x}}_{\c{E}^\ind{c},1}<\ls<
\stackrel{\circ}{\b{x}}_{\c{E}^\ind{c},k}<
\stackrel{\circ}{\b{z}}_{\c{E}^\ind{c}} $ that differ with one
coordinate only.
\end{proof}
Theorem \ref{t:main} allows to describe some properties of the set
of solutions
$$
\c{M}=\lfi\argmin_\b{x}P(\b{x})\rfi .
$$
\begin{cor}\label{c:struct}
Let $ \stackrel{\circ}{\b{x}}_{\c{E}^\ind{c}}\le
\stackrel{\circ}{\b{z}}_{\c{E}^\ind{c}} $ be any ordered pair of
vectors, then:

(i) for any $ \b{x}'_\c{D}\in
\c{M}({\stackrel{\circ}{\b{x}}_{\c{D}^\iind{c}}}) $ there exists a
solution $ \b{z}'_\c{D}\in
\c{M}({\stackrel{\circ}{\b{z}}_{\c{D}^\iind{c}}}) $ such that
$\b{x}'_\c{D}\le\b{z}'_\c{D}$.

(ii) For any $ \b{z}'_\c{D}\in
\c{M}({\stackrel{\circ}{\b{z}}_{\c{D}^\iind{c}}}) $ there exists a
solution $ \b{x}'_\c{D}\in
\c{M}({\stackrel{\circ}{\b{x}}_{\c{D}^\iind{c}}}) $ such that
$\b{x}'_\c{D}\le\b{z}'_\c{D}$.

(iii) The sets $ \c{M}({\stackrel{\circ}{\b{x}}_{\c{D}^\iind{c}}})
$ and $ \c{M}({\stackrel{\circ}{\b{z}}_{\c{D}^\iind{c}}}) $ have
minimal $ \underline{\b{x}}'_\c{D}, \underline{\b{z}}'_\c{D} $ and
maximal $ \overline{\b{x}}'_\c{D}, \overline{\b{z}}'_\c{D} $
elements.

(iv) The minimal and maximal elements are ordered, i.e. $
\underline{\b{x}}'_\c{D}\le \underline{\b{z}}'_\c{D} $ and $
\overline{\b{x}}'_\c{D}\le \overline{\b{z}}'_\c{D}. $
\end{cor}
The sentences (\emph{i}) and (\emph{ii}) follow directly from
Theorem \ref{t:main}. Sentence (\emph{iii}) is deduced from
Theorem \ref{t:main} for $ \stackrel{\circ}{\b{x}}_{\c{E}^\ind{c}}=
\stackrel{\circ}{\b{z}}_{\c{E}^\ind{c}}. $ And, at last,
(\emph{iv}) can be derived from (\emph{i,ii}) and definitions of
minimal and maximal elements.

\section{Modified SFM Algorithm}

The main idea of the \emph{Modified SFM algorithm} (MSFM) is to
determine (if possible) parts of a solution $\b{x}^*=\c{M}$ by use
separate parts $ P\{\c{D}_\ind{i}\}(\b{x}_{\c{D}_\iind{i}},
{\b{x}}_{\c{D}_\iind{i}^\iind{c}}), $ $
(\c{D}_\ind{i}\bigcup\c{D}_\ind{j}=\varnothing,\
\bigcap\c{D}_\ind{i}=\c{V}) $ of the original polynomial
$P(\b{x})$. For this purpose monotone dependence of local solutions
$\b{x}'_{\c{D}_\iind{i}}$ (in the sense of Corollary
\ref{c:struct}) on the frontier vector $
{\b{x}}_{\c{E}^\iind{c}_\iind{i}} $ can be exploited. If the vector
$ \b{x}'_{\c{D}_\iind{i}}\in \c{M}(\b{1}_{\c{E}^\ind{c}_\ind{i}}) $
has some coordinates $x'_j=0$, then there is a solution
$\b{x}^*\in\c{M}$ with the same coordinates $x^*_j=0$. Similarly,
if the solution $ \b{x}'_{\c{D}_\iind{i}}\in
\c{M}(\b{0}_{\c{E}^\ind{c}_\ind{i}}) $ is with coordinates
$x'_j=1$, then there is a solution $\b{x}^*\in\c{M}$ that has $
x^*_j=1. $

The strategy developed leads to success not for all Boolean
polynomials, but it is useful for polynomials that have local
dependence of variables. Such polynomials are quite often used in
applications, for instance, in Bayesian estimation or image
processing.

Describe the MSFM in details. Partition the set $\c{V}$ by
suitable sets $ \c{D}_\ind{i}(1),\ls,\c{D}_\ind{k_\iind{1}}(1) $ $
(\c{D}_\ind{i}(1)\bigcup\c{D}_\ind{j}(1)=\varnothing,\
\bigcap\c{D}_\ind{i}(1)=\c{V}). $ It follows from Proposition
\ref{p:pd} that $
\b{x}^*_{\c{D}_\iind{i}(1)}\in\c{M}(\b{x}^*_{\c{D}^\ind{c}_\iind{i}(1)}),
$ i.e. $\b{x}^*_{D_\iind{i}(1)}$ minimizes the polynomial $
P\{\c{D}_\ind{i}\}(\b{x}_{\c{D}_\iind{i}},
\b{x}^*_{\c{D}_\iind{i}(1)}). $ Corollary \ref{c:struct}
substantiates existence of triple of solutions
$$
\b{x}_{0,\c{D}_\iind{i}(1)}\in
\c{M}(\b{0}_{\c{D}^\ind{c}_\iind{i}(1)}),\quad
\b{x}_{1,\c{D}_\iind{i}(1)}\in
\c{M}(\b{1}_{\c{D}^\ind{c}_\iind{i}(1)})\quad \mbox{and}\quad
\b{x}'_{D_\iind{i}(1)}\in\c{M}(\b{x}^*_{\c{D}^\ind{c}_\iind{i}(1)})
$$
such that $ \b{x}_{0,\c{D}_\iind{i}(1)}\le
\b{x}'_{\c{D}_\iind{i}(1)}\le \b{x}_{1,\c{D}_\iind{i}(1)} $ (the
vectors can be found by the usual SFM algorithm). Therefore,
solutions $\b{x}'_{\c{D}_\iind{i}(1)}$ have coordinates
$x'_{\c{D}_\iind{i}(1),j}=0$ for  sets of indices $ \c{B}_i(1)=
\{j\in {\c{D}_\iind{i}(1)}\ |\ x_{1,{D_\iind{i}(1)},\ind{j}}=0\} $
as well as they have coordinates $x'_{D_\iind{i}(1),j}=1$ for sets
of indices $ \c{W}_i(1)= \{j\in {\c{D}_\iind{i}(1)}\ |\
x_{0,{\c{D}_\iind{i}(1)},\ind{j}}=1\}. $

Check the set $\c{R}(1)= \bigcup^{k_\ind{1}}_{i=1}
(\c{W}_\ind{i}(1)\bigcup \c{B}_\ind{i}(1))$. If it is empty we
either try another partition or start the standard SFM for the
original polynomial $ P(\b{x}). $

If $\c{R}(1)\neq\varnothing$, we identified the part of the
solution $ \b{x}^*_{\c{R}(1)}=\b{x}'_{\c{R}(1)} $. Set
$\c{V}(2)=\c{V}\bs \c{R}(1)$ and consider the reduced problem of
identification of
$$
\b{x}^*_{\c{V}(2)}= \argmin_{\b{x}_{\c{V}(2)}}
P(\b{x}_{\c{V}(2)},\b{x}^*_{\c{R}(1)}).
$$
Partition $ \c{V}(2)=\bigcup^{k_\ind{2}}_{i=1} \c{D}_\ind{i}(2),\
\c{D}_\ind{i}(2)\bigcap \c{D}_\ind{j}(2)=\varnothing $ and
estimate $\b{x}^*_{\c{D}_\ind{i}(2)}$ by solutions
$$
\b{x}_{0,\c{D}_\iind{i}(2)}\in
\c{M}(\b{0}_{\c{D}^\ind{c}_\iind{i}(2)\bs{\c{R}(1)}},
\b{x}^*_{{\c{R}(1)}})\quad \mbox{and}\quad
\b{x}_{1,\c{D}_\iind{i}(2)}\in
\c{M}(\b{1}_{\c{D}^\ind{c}_\iind{i}(2)\bs{\c{R}(1)}},
\b{x}^*_{\c{R}(1)}),
$$
using sets $ \c{B}_i(2)= \{j\in {\c{D}_\iind{i}(2)}\ |\
x_{1,{D_\iind{i}(2)},\ind{j}}=0\} $ and $ \c{W}_i(2)= \{j\in
{\c{D}_\iind{i}(2)}\ |\ x_{0,{\c{D}_\iind{i}(2)},\ind{j}}=1\} $
(Corollary \ref{c:struct} grounds existence of triple of vectors $
\b{x}_{0,\c{D}_\iind{i}(2)}\le \b{x}^*_{\c{D}_\ind{i}(2)}\le
\b{x}_{1,\c{D}_\iind{i}(2)} $ ). The algorithm is iterated at high
levels until the problem will be completely solved or until
$\c{R}(l)\neq\varnothing$. In the last case the standard SFM is
after applied. Sizes of partitioning sets $\c{D}_\iind{i}(l)$
should be sufficient to have $\c{R}(l)$ nonempty but not very
large to have small enough number of operations. It is unlikely
possible to estimate them in general case since for some
polynomials the MSFM is reduced to the usual SFM for any
partition. However, in many applied problems of image processing
or Bayesian estimation the finding suitable sizes of
$\c{D}_\iind{i}(l)$ is not difficult. In these problems energy
functions usually have lack of long distance interactions (or
phase transitions). It means solutions
$\b{x}^{*}_{\c{D}}$ for  the original image
and an image, which coincides with the original on sufficiently
large set $\c{D}\st\c{V}$ and differs from it on
$\c{D}^{\ind{c}}$, are approximately equal. By other words, for
real applications estimates of image details do not depend on its
surrounding objects. In this case the MSFM is preferable to the
standard SFM. The number operation required can be here even
$O(n)$.

The identification of $ \b{x}^*_{\c{B}_\iind{i}(l)\bigcup
\c{W}_\iind{i}(l)} $ can be done in concurrent mode. The number of
levels depend on the polynomial $P$. Usually, it is enough two or
three ones. Formal description of the MSFM is rather simple. It is
placed in Fig.1.

\nn \fbox {
\parbox{\textwidth}
{
\begin{tabbing}
\hspace{1.cm}\= \hspace{1cm}\= \textsf{MSFM}($P$):
\hspace{2cm}\=\hspace{1.5cm}\=\mbox{\hspace{6.7cm}}\\[6pt]
\mbox{\vspace{2cm}}\\
\> \textbf{Initialization:}\\[6pt]
\>\>$m$\> $\lla$\> $1\ $ the enumerator of levels\\[6pt]
\>\>$L$\> $\lla$\>a number of levels\\[6pt]
\>\>$\c{V}(1)$\> $\lla$\> $\c{V}$\\[6pt]
\>\>$\{\c{D}_\ind{i}(1)\}$\> $\lla$\> a partition of $\c{V}$\\[6pt]
\>\>$\{\b{x}_{0,\c{D}_\iind{i}(1)},\b{x}_{1,\c{D}_\iind{i}(1)}\}$\>
$\lla$\> such that
$\b{x}_{0,\c{D}_\iind{i}(1)}\le\b{x}_{1,\c{D}_\iind{i}(1)}$, by the SFM\\[6pt]
\>\>$\{\c{B}_\ind{i}(1),\c{W}_\ind{i}(1)\}$ \> $\lla$\>
coordinates of the part of the solution found\\[6pt]
\>\>$\c{R}(1)$\> $\lla$\>  $\bigcup^{k_\ind{1}}_{i=1}
(\c{W}_\ind{i}(1)
\bigcup \c{B}_\ind{i}(1))$\\[6pt]
\>\textbf{If} $\c{R}(1)=\varnothing$ \textbf{then} use the SFM, \emph{Stop}\\[6pt]
\>\>$x^*_{\c{D}_\iind{i}(1),j}$\> $\lla$\> $0\ $ for $j\in\c{B}_i(1)$\\[6pt]
\>\>$x^*_{\c{D}_\iind{i}(1),j}$\> $\lla$\> $1\ $ for $j\in\c{W}_i(1)$\\[6pt]
\>\textbf{While} $m< L$ \textbf{do}\\[6pt]
\>\>$m$ \> $\lla$\> $m+1$\\[6pt]
\>\>$\c{V}(m)$ \> $\lla$\> $\c{V}(m-1)\bs \c{R}(m-1)$\\[6pt]
\>\>$\{\c{D}_\ind{i}(m)\}$ \> $\lla$\> partition of $\c{V}(m)$\\[6pt]
\>\>$P_l(\b{x}_\ind{\c{V}(m)})$\> $\lla$\> the residual polynomial
$P(\b{x}_\ind{\c{V}(m)},\b{x}^*_\ind{\c{V}^\iind{c}(m)})$\\
\>\>$\{\b{x}_{0,\c{D}_\iind{i}(m)},\b{x}_{1,\c{D}_\iind{i}(m)}\}$\>
$\lla$\> such that
$\b{x}_{0,\c{D}_\iind{i}(m)}\le\b{x}_{1,\c{D}_\iind{i}(m)}$, by the SFM\\[6pt]
\>\>\>\> for $P_l(\b{x}_\ind{\c{V}(m)})$\\[6pt]
\>\>$\{\c{B}_\ind{i}(m),\c{W}_\ind{i}(m)\}$ \> $\lla$\>
coordinates of the part of the solution found\\[6pt]
\>\>\hspace{-0.5cm}\textbf{If} $\bigcup^{k_\ind{m}}_{i=1}
\lb\c{W}_\ind{i}(m) \bigcup \c{B}_\ind{i}(m)\rb=\varnothing$
\textbf{then} use the SFM, \emph{Stop}\\[6pt]
\>\>$\c{R}(m)$\> $\lla$\>  $\c{R}(m-1)\bigcup
\lb\bigcup^{k_\ind{m}}_{i=1} \lb\c{W}_\ind{i}(m)
\bigcup \c{B}_\ind{i}(m)\rb$\rb\\[6pt]
\>\>$x^*_{\c{D}_\iind{i}(m),j}$\> $\lla$\> $0\ $ for $j\in\c{B}_i(m)$\\[6pt]
\>\>$x^*_{\c{D}_\iind{i}(m),j}$\> $\lla$\> $1\ $ for $j\in\c{W}_i(m)$\\[6pt]
\>\textbf{If} $\c{R}(L)\neq\varnothing$
\textbf{then} use the SFM\\[6pt]
\>\textbf{Return} $\b{x}^*$\\[6pt]
\>\textbf{End}
\end{tabbing}
} } \vspace{0.5cm}
\begin{center}
{ Fig.1 }
\end{center}

\newpage
\section{Representability of Some Submodular Polynomials by Graphs }

Recently  Kolmogorov\&Zabih \cite{Kol01} have described some class
of functions $V(\b{x})$ in Boolean variables
$\b{x}=(x_1,\ldots,x_n)$ that permit minimization by the graph cut
technique. For convenience we reformulate and then use their
results in terms of Boolean polynomials. We also represent some
submodular Boolean polynomials by graphs to make possible their
efficient minimization via graph gut technique (recall that
finding the minimum graph cut requires $O(N^3)$ operations
on number $N$ of nodes).

Let the network ${\Gamma}=(\c{S},\c{A})$ consist of $N+2$ numbered
nodes ${\c{S}}=\{0,\ls,N+1\}$, where $s=0$ is the \emph{source},
$t=N+1$ is the \emph{sink} and ${\c{V}^{\, \prime}}=\{1,\ls,N\}$ are usual
nodes. The set of directed arcs is $\c{A}=\{(i,j)\ :\ i,j\in
{\c{S}}\}$. The capacities of arcs are denoted by $d_{i,j}>0$.

The \emph{cut} (we use this term as abbreviation of the term
\emph{s-t cut}) is a partition of the set ${\c{S}}$ by two
disjoint sets $W,B$ such that $s\in W$ and $t\in B$. For any cut
$(W,B)$ nodes $i\in W$ will be labeled by $z_i=1$ and nodes $i\in
B$ -- by $z_i=0$ so that there is one-to-one map between cuts and
Boolean vectors $\b{z}=(z_1,\ldots,z_N)$.

\subsection{Graph representability}
In \cite{Kol01} the following definition of graph representability
is given.
\begin{f}\label{d:kzr}
A function $V(\b{x})$ of $n$ Boolean variables is called graph
representable if there exists a network $\Gamma=(\c{S},\c{A})$
with number of nodes $N\ge n$ and the cost function $C(\b{z})$
such that for some constant \rf{c}
$$
V(\b{x})=\min_\ind{z_\iind{n+1},\ls,z_N} C(x_1,\ls,x_n,z_{n+1},\ls
,z_N) +\rf{c}.
$$
\end{f}
Costs of the graph cuts $C(\b{z})$ are equal to values of the
quadratic Boolean polynomial $p_{}(\b{z})$ of the form (see,
\cite{GPS89,PR75}).
\begin{equation}\label{e:pg}
C(\b{z})=p_{}(\b{z})= \sum_{i\in \c{V}}d_{s,i}(1-z_i)+ \sum_{i\in
\c{V}}d_{i,t}z_i+ \sum_{i,j\in \c{V}}d_{i,j}(z_i-z_j)z_i +\rf{c},
\quad d_{i,j},\rf{c}\ge 0.
\end{equation}
Therefore, we can give an equivalent definition of the graph
representability.
\begin{f}\label{d:gr}
We say a function $V(\b{x})$ of $n$ Boolean variables can be
represented by graph if there exists a quadratic polynomial
$\b{p}(\b{z})$ of $N\   (N\ge n)$ Boolean variables of the form
(\ref{e:pg}), which satisfies the equality
$$
V(\b{x})=\min_\ind{z_\iind{n+1},\ls,z_N}
\b{p}(x_1,\ls,x_n,z_{n+1},\ls ,z_N)+\rf{c}.
$$
\end{f}
 The class $\c{F}^2$ of
Boolean functions that are often used in applications and that
are, actually, quadratic Boolean polynomials
\begin{equation}\label{e:F2}
V\in\c{F}^2 \Leftrightarrow V(\b{x})=P_V(\b{x})= \sum_\ind{1\le l
\le n}a_\ind{l}x_{l}+\ \ \sum_\ind{1\le
{l}_\iind{1}<{l}_\iind{2}\le n}\ \
 a_\ind{l_\iind{1},l_\iind{2}}x_\ind{l_\iind{1}}x_\ind{l_\iind{2}}
+\rf{c}.
\end{equation}
was considered in \cite{Kol01}. The characterization of the graph
representability of $V\in\c{F}^2$ in terms of Boolean polynomials
has been  done. It can be used, in particular, to minimize
functions of integer variables.
\begin{thm}\label{t:rep}
A function $V\in\c{F}^2$ is graph representable if and only if all
quadratic coefficients of its polynomial representation
$P_V(\b{x})$ are nonpositive, i.e. $
a_\ind{l_\iind{1},l_\iind{2}}\le 0$. Therefore, in order to
function $V\in\c{F}^2$ be graph representable it is necessary an
sufficient it will be submodular.
\end{thm}
The sufficiency of the Theorem follows immediately from Definition
\ref{d:gr} and Proposition~\ref{p:prop}~a). The necessity will be
proven below for more general case. The following result has been
proven in \cite{Kol01} by graph cut technique.
\begin{thm}\label{t:s-nes}
In order to  polynomial $P(\b{x})$ of degree $m$ can be
represented by a graph  it should be submodular.
\end{thm}
\begin{proof}
To prove this sentence we use representation (\ref{e:pij}) of the
polynomial $P(\b{x})$ at page~\pageref{pij}. For arbitrary couple of
indices $i,j$ let us fixed the vector $\b{x}_{\{i,j\}^\iind{c}}$.
Since $P(\b{x})$ is graph representable the quadratic polynomial $
x_ix_jP_{i,j}(\b{x}_{\{i,j\}^\iind{c}}) $ of two variables
$x_i,x_j$ is graph representable, as well. That is
$$
x_ix_jP_{i,j}(\b{x}_{\{i,j\}^\iind{c}})=
\min_\ind{z_\iind{n+1},\ls,z_N} \b{p}(x_1,x_2,z_{n+1},\ls
,z_N)+\rf{c}
$$
(in general, the quadratic submodular polynomial $\b{p}$ and
$\rf{c}$ depend on $\b{x}_{\{i,j\}^\iind{c}}$ but it is fixed now).
Suppose $P_{i,j}(\b{x}_{\{i,j\}^\iind{c}}) > 0$. Then $\b{p}$ has
at least two values of the argument $(x_1,x_2,z_{n+1},\ls ,z_N)$
that minimize the polynomial globally. One of them is with
coordinates $(0,1,\ls)$, and the other is of the form $(1,0,\ls)$.
Because of submodularity of $\b{p}$ it follows from Corollary
\ref{c:struct}  (i)  there exists the third vector, which
minimizes $p$ globally. This vector has coordinates $(1,1,\ls)$.
The existence of such the vector contradicts with the assumption
$P_{i,j}(\b{x}_{\{i,j\}^\iind{c}}) > 0$. Therefore,
$P_{i,j}(\b{x}_{\{i,j\}^\iind{c}}) \le 0$ for any couple of
indices $i,j$ and  vector $\b{x}_{\{i,j\}^\iind{c}}$ and the
sentence of the Theorem follows from Corollary \ref{p:prop} a).
\end{proof}

\subsection{Graph represtability of some subset of submodular
polynomials}

The submodularity is necessary condition in order to a Boolean
polynomial can be represented by  a graph. Below we prove
representability by graphs of some subset $\c{P}_{suf}$
of submodular
polynomials. To define this subset we use formula
(\ref{e:pij}) at page \pageref{pij}.
Let
$
b_{l_\iind{1},\ls,l_\iind{m}}(i,j)=
a_{l_\iind{1},\ls ,i,\ls ,j,\ls ,l_\iind{m}}
$,
$
(l_\iind{1}<\ls< i<\ls <j<\ls<l_\iind{m},
\
m=0\div n-2)$
be  coefficients of
$P_{i,j}(\b{x}_{\{i,j\}^\iind{c}})$
and $b^{+}_{l_\iind{1},l_\iind{2},\ls,l_\iind{m}}(i,j)$
be the positive ones (note that $b(i,j)=a_\ind{i,j}<~0$).
Let the set of polynomials
$$
\c{P}_{suf}=\lfi P\ \ {\Big |}\ \
-a_{i,j}\ +\ \sum_{m=1}^{n-2}\
\sum_{l_\iind{1}<l_\iind{2}<\ls<l_\iind{m}}
b^{+}_{l_\iind{1},l_\iind{2},\ls,l_\iind{m}}(i,j)
\le 0,\ \ \forall i,j\in\c{V}\rfi.
$$
This set is a proper subset of all
submodular polynomials $\c{P}_{submod}$
(see Proposition~\ref{p:prop}~a)).
It  contains, for instance,
the set $ \c{F}^m_{-}= \{ P\ |\
a_{l_\iind{1},l_\iind{2},\ls,l_\iind{k}}\le 0,\ k=2\div m \} $
of Boolean polynomials with all coefficients of nonlinear
monomials nonpositive and the set
$ \c{F}^m_{+}= \{ P\ |\ P\in\c{P}_{submod},\
a_{l_\iind{1},l_\iind{2},\ls,l_\iind{k}}\ge 0,\ k=3\div m \} $
of submodular polynomials with all coefficients of order more
than two nonnegative.
\begin{rem}
Note the condition
$$
-a_{i,j}\ +\ \sum_{m=1}^{n-2}\
\sum_{l_\iind{1}<l_\iind{2}<\ls<l_\iind{m}}
b^{+}_{l_\iind{1},l_\iind{2},\ls,l_\iind{m}}(i,j)
\le 0,\ \ \forall i,j\in\c{V}
$$
is necessary and sufficient in order to
the polynomial, which has all coefficients of order
more than two positive (i.e. $P\in\c{P}_{submod}$), be submodular.
\end{rem}
Polynomials of those types can be used, for
instance, in image processing for recognition of images or
Bayesian estimation. To prove representability of
$P\in\c{P}_{suf}$ by graphs we need two lemmas,
which give a constructive tool to build a corresponding graph as
well.
\begin{lem}\label{l:neg}
For any natural $m$ and real $a<0$ the Boolean polynomial
$ax_1x_2\ls x_m$ is graph representable.
\end{lem}
\nn \fbox{\includegraphics[width=4.7cm,height=5.3cm] {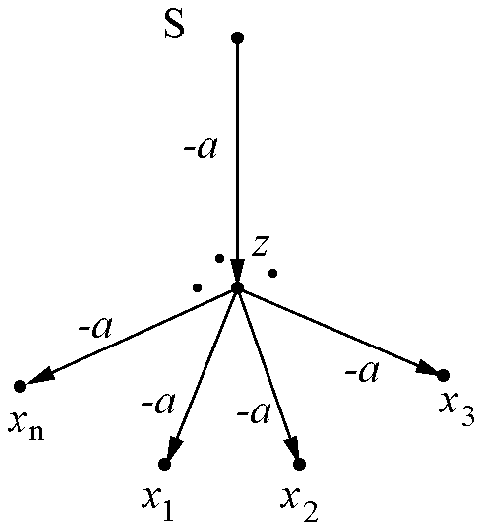}}

\vspace{-5.6cm} \hfill\parbox{10cm} { The graph that represents
$ax_1x_2\ls x_m$ is drawn in Fig.2. It needs one additional
node~$z$. The corresponding polynomial is
$p(x_1,x_2,\ls,x_m,z)=-a\sum_{i=1}^m (z-x_i)z +az$. If at least
one $x_i=0$, then $\min_{z}p(x_1,x_2,\ls,x_m,z)=0$ otherwise
$\min_{z}p(1,1,\ls,1,z)=a$. Hence, $ax_1x_2\ls x_m$
$=\min_{z}p(x_1,x_2,\ls,x_m,z).$

The monomials of order greater than $2$ with positive coeffisients
can not be represented  by graphs directly. To do it we should
accompany those positive monomials by quadratic items with
large enough negative coefficients.
It requires additional nodes in  graphs that represent
the monomials.
}

\nn \hspace{2cm} Fig.2


\begin{lem}\label{l:pos}
For any natural $m>2$ and $a>0$ the Boolean polynomial
$P(x_1,x_2,\ls, x_n)=ax_1x_2\ls x_m-a\sum_{1\le i<j\le m}x_i x_j$
is graph representable.
\end{lem}

\begin{proof}
In order to prove the Lemma we construct for each $m>2$ a
quadratic Boolean polynomial that specifies graph (the examples of
graphs  for $m=3\div 6$ are drawn in Fig.3 at page \pageref{p:pl}).
The number of
additional variables of this polynomial is $l=\dfrac{m-1}{2}$ if
$m$ is odd and is $l=\dfrac{m-2}{2}$ if $m$ is even. It is of the
form
\begin{equation}\label{e:pfm}
p(x_1,\ls,x_m,z_1,\ls,z_l)=\sum_{j=1}^l
b_j\sum_{i=1}^m(z_j-x_i)z_j +\sum_{j=1}^l e_jz_j
\end{equation}
with all $b_j>0$ and $e_j<0$, such that the minimum of $p$ with
respect to  $z_1,\ls, z_l$ is equal to corresponding  values of
$P$. Those values are
$$
P(x_1,x_2,\ls, x_m)=
\begin{cases}
0,              &\mbox{if all}\ x_i=0,\\
0,              &\mbox{if}\ x_\ind{i_\iind{1}}=1\
\mbox{and}\ x_\ind{i_\iind{2}}=\ls= x_\ind{i_\iind{m}}=0\\
-a\dbinom{k}{2},&\mbox{if}\ x_\ind{i_\iind{1}}=\ls=
x_\ind{i_\iind{k}}=1\
\mbox{and}\ x_\ind{i_\iind{k+1}}=\ls= x_\ind{i_\iind{m}}=0\\
\\
-a\dbinom{m}{2}+a,&\mbox{if}\ x_1=\ls= x_m=1.
\end{cases}
$$
Let the vector $\b{1}_k =
(\overbrace{\underbrace{1,\ls,1}_{k},0,\ls,0}^m)$. The polynomial
$p$ (as well as $P$) is symmetrical with respect to $x_1,\ls,x_m$.
It has the same value $p(\b{1}_k)= \sum_{j=1}^l( (m-k)b_j+e_j)z_j$
for all permutations of $\b{1}_k$. Therefore, the set of
inequalities should hold true
\begin{equation}\label{e:ine}
\begin{cases}
\sum\limits_{j=1}^l ( mb_j+e_j)z_j\ & \ge 0\\
\sum\limits_{j=1}^l ( (m-1)b_j+e_j)z_j\ & \ge 0\\
\sum\limits_{j=1}^l ( (m-2)b_j+e_j)z_j\ & \ge -a\\
\vdots\\
\sum\limits_{j=1}^l ( (m-k)b_j+e_j)z_j\ & \ge -a\dbinom{k}{2}\\
\vdots\\
\sum\limits_{j=1}^l (b_j+e_j)z_j\ & \ge -a\dbinom{m-1}{2}\\
\sum\limits_{j=1}^l e_jz_j\ & \ge -a\dbinom{m}{2}+a
\end{cases}
\end{equation}
together with the additional condition: \vspace{0.5cm}

(M) \hspace{0.5cm} \emph{There exist vectors $\{\b{z}^*_k\}$ that
turn corresponding inequalities into equality.}

\vspace{0.5cm}

For any odd $m>2$ we propose one of possible collections of
variables $b_i,e_i$ that satisfy (\ref{e:ine}) for Boolean
variables $(z_1,\ls,z_m)$ as well as condition (M). The collection
is identified by the following system of equations
\begin{alignat}{1}\label{e:coefo}
\begin{cases}
\big( (m-2)b_1+e_1\big)\ & = -a  \dbinom{2}{2}\\
\big( (m-3)b_1+e_1\big)\ & = -a  \dbinom{3}{2}\\
\big( (m-4)b_1+e_1\big)+
\big( (m-4)b_2+e_2\big)\ & = -a \dbinom{4}{2}\\
\big( (m-5)b_1+e_1\big)+
\big( (m-5)b_2+e_2\big)\ & = -a \dbinom{5}{2}\\
\vdots\\
\sum\limits_{j=1}^{k} \big( (m-2k)b_j+e_j\big)\ & = -a \dbinom{2k}{2}\\
\sum\limits_{j=1}^{k} \big( (m-2k-1)b_j+e_j\big)\ & = -a \dbinom{2k+1}{2}\\
\vdots\\
\sum\limits_{j=1}^l \big(b_j+e_j\big)\ & = -a \dbinom{m-1}{2}\\
\sum\limits_{j=1}^l e_j\ & = -a \dbinom{m}{2}+a
\end{cases}
\end{alignat}
It is
\begin{align}\label{e:solo}
b_1&=b_2=\ls=b_{l-1}=2a,\qquad b_l=a, \nonumber\\
\\
e_j&=-(2m-4j+1)a,\ (j=1\div l-1),\qquad e_l=-2a\nonumber.
\end{align}
Since  the inequalities
\begin{align*}
&(m-1)b_j+e_j>0,\quad j=1\div l,\\
&(m-2k)b_j+e_j<0,\quad (m-2k-1)b_j+e_j<0,\quad \mbox{if}\quad
j\le k,\ k=1\div l-1,\\
&(m-2k)b_j+e_j>0,\quad (m-2k-1)b_j+e_j>0,\quad \mbox{if}\quad
j>k,\ k=1\div l-1,\\
&b_j+e_j<0,\ j=1\div l,\\
&e_j<0,\quad j=1\div l.
\end{align*}
hold true, the solution (\ref{e:solo}) actually satisfies
inequality (\ref{e:ine}) together with condition (M) for  vectors
$\b{z}^*_{1}=\b{z}^*_{2}=\b{0},\ \
\b{z}^*_{3}=\b{z}^*_{4}=\b{1}_{1},\ \ \ls\ \ ,
\b{z}^*_{2k+1}=\b{z}^*_{2k+2}=\b{1}_{k},\ \ls,\ \
\b{z}^*_{m}=\b{z}^*_\ind{m+1}=\b{1}_{l}$ .

For even $m>2$ the collection  of $l=\dfrac{m-2}{2}$ variables
$b_i,e_i$ can be found as a solution of the system of equation
\begin{alignat}{1}\label{e:coefe}
\begin{cases}
\big( (m-2)b_1+e_1\big)\ & = -a  \dbinom{2}{2}\\
\big( (m-3)b_1+e_1\big)\ & = -a  \dbinom{3}{2}\\
\big( (m-4)b_1+e_1\big)+
\big( (m-4)b_2+e_2\big)\ & = -a \dbinom{4}{2}\\
\big( (m-5)b_1+e_1\big)+
\big( (m-5)b_2+e_2\big)\ & = -a \dbinom{5}{2}\\
\vdots\\
\sum\limits_{j=1}^{k} \big( (m-2k)b_j+e_j\big)\ & = -a \dbinom{2k}{2}\\
\sum\limits_{j=1}^{k} \big( (m-2k-1)b_j+e_j\big)\ & = -a \dbinom{2k+1}{2}\\
\vdots\\
\sum\limits_{j=1}^l \big(2b_j+e_j\big)\ & = -a \dbinom{m-2}{2}\\
\sum\limits_{j=1}^l \big(b_j+e_j\big)\ & = -a \dbinom{m-1}{2}
\end{cases}
\end{alignat}
that looks similar to (\ref{e:coefo}) but consists of only $m$
equations, since the equation
$$
\sum\limits_{j=1}^l e_jz_j\  = -a\dbinom{m}{2}+a
$$
follows from the last two ones of system (\ref{e:coefe}). The
solution of the system is
\begin{equation}\label{e:sole}
b_{j}=2a,\qquad e_j=-(2m-4j+1)a,\qquad j=1\div l.
\end{equation}
Again the inequalities
\begin{align*}
&(m-1)b_j+e_j>0,\quad j=1\div l,\\
&(m-2k)b_j+e_j<0,\quad (m-2k-1)b_j+e_j<0,\quad \mbox{if}\quad
j\le k,\ k=1\div l,\\
&(m-2k)b_j+e_j>0,\quad (m-2k-1)b_j+e_j>0,\quad \mbox{if}\quad
j>k,\ k=1\div l,\\
\end{align*}
hold true and, therefore, (\ref{e:ine}) together with condition (M)
is fulfilled for  vectors $\b{z}^*_{1}=\b{z}^*_{2}=\b{0},\ \
\b{z}^*_{3}=\b{z}^*_{4}=\b{1}_1,\ \ \ls\ \ ,
\b{z}^*_{2k+1}=\b{z}^*_{2k+2}=\b{1}_k,\ \ls,\ \
\b{z}^*_{m-1}=\b{z}^*_{m}=\b{z}^*_{m+1}=\b{1}_l$ .

So, the coefficients of the polynomial $\b{p}$ (see (\ref{e:pfm}))
that represents $P$ for odd and even $m$ are determined by
(\ref{e:solo}) and (\ref{e:sole}).
\end{proof}
Lemmas  \ref{l:neg} and \ref{l:pos} allow to prove
each
$P\in\c{P}_{suf}$ can be represented by a graph.

\begin{center}
\fbox{\includegraphics[width=6cm,height=7.2cm] {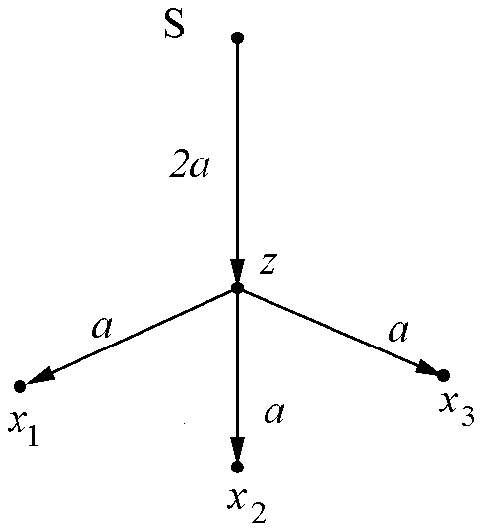}}
\hspace{2cm}
\fbox{\includegraphics[width=6cm,height=7.2cm] {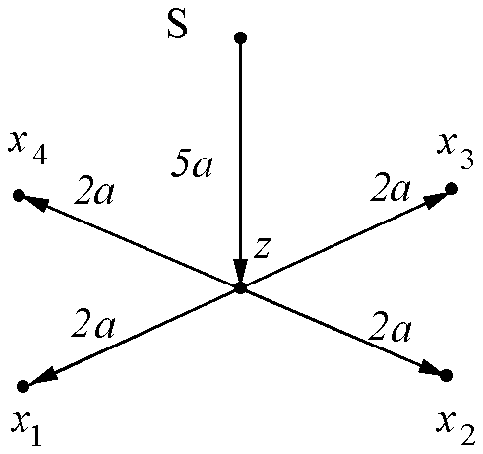}}
\end{center}

\begin{center}
\fbox{\includegraphics[width=6cm,height=7.2cm] {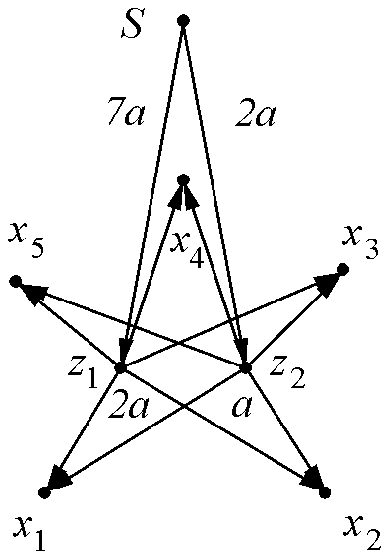}}
\hspace{2cm}
\fbox{\includegraphics[width=6cm,height=7.2cm] {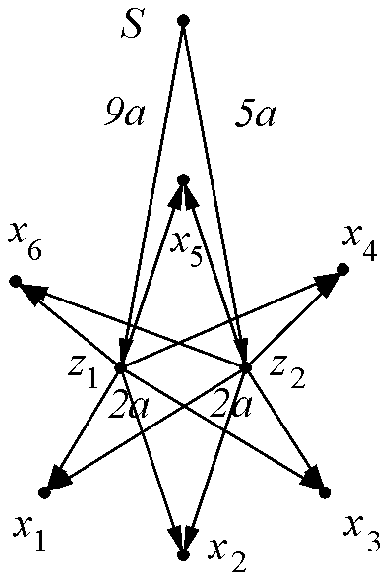}}
\end{center}

\begin{center}
Fig.3
\end{center}
\label{p:pl}

\begin{thm}
The polynomials $P\in\c{P}_{suf}$ are graph
representable.
\end{thm}

In order to proof the Theorem it is enough to use Lemma \ref{l:neg}
to represent nonlinear items of $P$ with negative coefficients and,
then, to exploit Lemma \ref{l:pos} to determine that all items
with positive coefficients can be accompanied by quadratic
monomials and after be represented by graphs as well.

Let the Boolean polynomial $P(x_\ind{1},\ls,x_\ind{n})\in\c{P}_{suf}$
of degree $m$ consists of $L_{h}$ monomials of degree
more than two, then it can be represented by a graph
with  $\tfrac{L_{h}(m-1)}{2}+n$ nodes.
The worst case is when $P(x_\ind{1},\ls,x_\ind{n})\in\c{F}^m_{+}$.
The polynomial
$P(x_\ind{1},\ls,x_\ind{n})\in\c{F}^m_{-}$ can be
represented by a graph with $L_{h}+n$ nodes.
Hence, the global minimization of
$
P(x_\ind{1},\ls,x_\ind{n})\in\c{P}_{suf}
$
require less than $O\lb\tfrac{L_{h}(m-1)}{2}+n\rb^3$
operations, though for some classes of polynomials
that occur in applications the number of operations can be
reduced even to $O\lb\tfrac{L_{h}(m-1)}{2}+n\rb$ (see, \cite{Z02}).

\end{document}